\documentclass[12pt]{iopart}
 \usepackage{amssymb}
 \usepackage{epsfig}
  \usepackage{pict2e}
  \usepackage{graphicx}
 \newtheorem{Lemma}{Lemma}
 \newtheorem{Proposition}[Lemma]{Proposition}

\newcommand{\EP}{empire process}
 
  \newcommand{\Dcrit}{{D}_{\mbox{{\tiny crit}}}}
  \newcommand{\gamcrit}{{\gamma}_{\mbox{{\tiny crit}}}}
 \newcommand{\area}{\mathrm{area}}
 \newcommand{\peri}{\mathrm{peri}}
  \newcommand{\sfrac}[2]{{\textstyle\frac{#1}{#2}}}
  \newcommand{\RR}{\mbox{${\mathcal R}$}}
  \newcommand{\CC}{\mbox{${\mathcal C}$}}
  \newcommand{\FF}{\mbox{${\mathcal F}$}} 

\begin{document}
 \title[Empires and percolation]{Empires and percolation: stochastic merging of adjacent regions}
 \author{D J Aldous\footnote{Research supported by
 N.S.F Grant DMS-0704159}, J R Ong  and W Zhou}
 \address{Department of Statistics\\
 367 Evans Hall \#\  3860\\
 U.C. Berkeley CA 94720}
 \ead{ aldous@stat.berkeley.edu
 \\ www.stat.berkeley.edu/users/aldous}

\begin{abstract}
We introduce a stochastic model in which adjacent planar regions $A, B$ merge stochastically at some  rate $\lambda(A,B)$, and observe analogies with the well-studied topics of mean-field coagulation and of bond percolation.   Do infinite regions appear in finite time?  We give a simple condition on $\lambda$ for this {\em hegemony} property to hold, and another simple condition for it to not hold, but there is a large gap between these conditions, which includes the case $\lambda(A,B) \equiv 1$.  
For this case, a non-rigorous analytic argument and simulations suggest hegemony.  
\end{abstract}



 \vspace{0.4in}

 
  \maketitle

  \newpage

 \section{Introduction}
 \label{sec-INT}
We study random processes defined as follows. \\
(i) at each time $t$ the plane is partitioned into polygonal regions\\  (ii) as $t$
increases, adjacent regions $A, B$ merge into one region  $A \cup B$ stochastically at some
rate (probability per unit time) 
$\lambda(A,B)$.

\vspace{0.1in}
\noindent
Here $\lambda(A,B)$ is a specified function of the geometry of regions $A$ and $B$,
typically a simple formula involving quantities such as the areas
$\area(A), \area(B)$, the perimeters (boundary lengths) 
$\peri(A), \peri(B)$ and the length $L(A,B)$ of the boundary between $A$ and $B$. 
The {\em adjacency} condition is enforced by the assumptions 
\[ \mbox{if $L(A,B) > 0$ then $\lambda(A,B) > 0$; \ \ \ if  $L(A,B) = 0$ then
$\lambda(A,B) = 0$. }\]
This is intended as an abstract model for spatial growth via merging -- we mentally
picture countries merging into empires.
Despite its conceptual simplicity, this model has apparently never been studied, so
let us name it the {\em \EP}. 
The purpose of this article is to initiate its study by setting out some simple
analytic and simulation results, describing connections with other models, 
and providing suggestions for further study and methodology.

\subsection{Analogy with mean-field coalescence}
Consider moving particles of different masses  $x$ 
(we take the continuous case $0<x<\infty$), where particles of masses $x$ and $y$
may coalesce when they meet. 
Ignoring geometry (positions and motion of particles) and simply assuming  there is
a "relative propensity to meet and coalesce" function $K(x,y)$ which applies to all
pairs of particles, leads one to formulate a  ``mean-field" model.  
Write $f(x,t) dx$ for the density per unit volume at time $t$ of particles with mass
in $[x,x +dx]$.  
The density function $f(x,t)$ satisfies the well-known
{\em Smoluchowski coagulation equation} \cite{drake,me78,norris99}  \begin{equation}
\sfrac{d}{dt} f(x,t) 
= \sfrac{1}{2} \int_0^x \! K(y,x-y) f(y,t) f(x-y,t) dy
 -  f(x,t) \int_0^\infty  \! K(x,y) f(y,t) dy .
\label{smol2}
\end{equation}
In the \EP\  the conserved quantity is area, so there is an analogous density
$f(a,t)$ of area-$a$ regions per unit area.
But -- even when the merger rate function is of the form $\lambda(A,B) = K(\area(A),
\area(B))$ -- 
we cannot write down autonomous equations analogous to (\ref{smol2}) for
$\sfrac{d}{dt} f(a,t)$, because spatial relationships, and the stochastic nature of
mergers, matter in our model.

Versions (now called {\em stochastic  coalescents}) with finite total mass and
stochastic merging 
but no geometry have been studied  \cite{me78, MR2253162}.
And there is recent progress \cite{MR2350433} in rigorous verification of the
underlying presumption that in models of spatial diffusion and merging, the limiting
(low density of massive particles) behaviour is as predicted by the Smoluchowski
coagulation equation.
Our \EP\ model is just different, but our general theme of relating qualitative
properties of the process to properties of the rate function  is obviously parallel
to the same theme in mean-field coalescence.

\subsection{Analogy with bond percolation}
To each edge $e$ of the square lattice assign a random time $T_e$, with
Exponential(1) distribution 
$P(T_e \leq t) = 1-e^{-t}$,  at which the edge becomes ``open".    The configuration
of open edges at time $t$ is just the usual bond percolation \cite{gri99,MR2283880} 
process with 
$p =  1-e^{-t}$.  
There is an associated spanning forest process in which is included only those edges
which (upon becoming open) link two distinct open components.  
A connected component of the bond percolation process with $m+1$ vertices 
contains a tree (within the spanning forest process) with $m$ edges, 
and this connected component can be identified, as in Figure 1, with a
region of area $m$ in the dual square grid.  

\setlength{\unitlength}{0.15in}
\begin{picture}(10,9)(-8,-2)
\put(0,2){\circle*{0.5}}
\put(0,4){\circle*{0.5}}
\put(2,0){\circle*{0.5}}
\put(2,2){\circle*{0.5}}
\put(2,4){\circle*{0.5}}
\put(4,0){\circle*{0.5}}
\put(4,2){\circle*{0.5}}
\put(4,4){\circle*{0.5}}
\put(6,2){\circle*{0.5}}
\put(8,2){\circle*{0.5}}
\put(10,2){\circle*{0.5}}
\put(8,0){\circle*{0.5}}
\thinlines
\put(-1,1){\line(0,1){4}}
\put(-1,5){\line(1,0){6}}
\put(5,5){\line(0,-1){2}}
\put(5,3){\line(1,0){6}}
\put(-1,1){\line(1,0){2}}
\put(1,1){\line(0,-1){2}}
\put(1,-1){\line(1,0){4}}
\put(5,-1){\line(0,1){2}}
\put(5,1){\line(1,0){2}}
\put(7,1){\line(0,-1){2}}
\put(7,-1){\line(1,0){2}}
\put(9,-1){\line(0,1){2}}
\put(9,1){\line(1,0){2}}
\put(11,1){\line(0,1){2}}
\linethickness{0.5mm}
\put(0,2){\line(0,1){2}}
\put(0,2){\line(1,0){2}}
\put(2,0){\line(0,1){4}}
\put(4,0){\line(0,1){4}}
\put(2,0){\line(1,0){2}}
\put(4,2){\line(1,0){6}}
\put(8,2){\line(0,-1){2}}
\end{picture}

{\bf Figure 1.}  {\small The region associated with a component in bond percolation.}

\vspace{0.12in}
\noindent
The number of edges in the lattice between two adjacent components of the bond
percolation process equals the boundary length between the two associated dual
regions.  
In the bond percolation process as $t$ increases to $t + dt$, each closed edge has
chance $dt$ to become open, so the bond percolation components merge at rate 
``number of closed edges linking the components", and so 
in the dual process two regions $A,B$ merge at rate  $L(A,B)$.  
So we have shown
\begin{quote}
In the \EP\ with $\lambda(A,B) = L(A,B)$, started at $t = 0$ with the unit squares
of the square grid, the areas of regions at time $t$ are
distributed as the numbers of edges of the connected components of bond percolation
with $p(t) =   1-e^{-t}$.
\end{quote}
The celebrated result  \cite{sykes,kesten} that the critical value in bond
percolation equals $1/2$ implies that, in the particular \EP\  above, infinite
regions appear at time $\log 2$.

\section{Hegemonic or not?}
\label{sec-heg}
Whether or not the qualitative property
\begin{equation}
\mbox{infinite-area regions appear in finite time} 
\label{def-heg}
\end{equation}
holds is perhaps the most interesting question to ask about \EP es. In percolation
theory the analogous property is described by phrases like ``supercritical" or
``percolation occurs", while in coagulation theory it is called ``gelation", but to
maintain the visualization of  empires let us call (\ref{def-heg}) the {\em
hegemonic} case.   
(See Notes below for details of interpretation of (\ref{def-heg})). Can we give
conditions, in terms of the merger rate function
$\lambda(A,B)$, for whether the \EP\ is hegemonic or non-hegemonic?  Here is one
easy result.
\begin{Proposition}
\label{P1}
(i) Suppose (for some $c>0$) we have
\[ \lambda(A,B) \geq c \ L(A,B) \mbox{ for all } A, B . \]
Then the \EP\ is hegemonic. \\
(ii) Suppose (for some $C<\infty$) we have 
\[ \sum_j \lambda(A,B_j)  \peri(B_j) \leq C\  \peri(A)  \mbox{ for all $A$ with
adjacent regions $(B_j)$}. \]
Then the \EP\ is not hegemonic. 
\end{Proposition}
See (\ref{ex-explicit}) for an explicit example of (ii).

Before giving the simple proofs,  we should point out a subtlety whose resolution
will be useful. 
 Consider two \EP es, where the second has faster rates than the first:
\begin{equation}
  \lambda_1(A,B) \leq \lambda_2(A,B)\mbox{ for all } A, B . 
  \label{rel-0}
  \end{equation}
The assertion
\begin{eqnarray}
\mbox{if the $\lambda_1$ \EP\ is hegemonic and (\ref{rel-0}) holds then}\nonumber \\
\mbox{the $\lambda_2$ \EP\ is hegemonic} \label{heg12}
\end{eqnarray}
seems plausible at first sight, but 
 we suspect it is not always true (see Notes below). 
 What is true?  
 If we can couple (define simultaneously) the two processes in such a way
that the natural 
 refinement-coarsening relationship is maintained 
 (each region of the first is a subset of a region of the second) then the
conclusion 
 ``if the first process is hegemonic then the second process is also 
hegemonic" is certainly correct. 
 And the natural condition which ensures the relation can be maintained is
the following.
 \begin{equation}
 \sum_{i,j} \lambda_1(A_i,B_j) \leq \lambda_2(A,B), \quad 
 \mbox{ for all partitions $(A_i)$ of $A$ and $(B_j)$ of $B$}. 
 \label{rel-12}
 \end{equation}
 This condition includes the case $A_1 = A, B_1 = B$ and so (\ref{rel-12})
is a
stronger assumption than (\ref{rel-0}).

To investigate when condition (\ref{rel-12}) might hold, consider a rate function
with the {\em superadditive} property
\begin{equation}
\lambda_1(A_1,B) + \lambda_1(A_2,B)\leq \lambda_1(A_1 \cup A_2,B), \quad 
 \mbox{ for all $B$ and disjoint $(A_1, A_2)$}. 
 \label{rel-3}
 \end{equation}
This condition implies that (\ref{rel-12}) holds with $\lambda_2$ replaced by $
\lambda_1$.  So if we assume both (\ref{rel-0}) and superadditivity of $\lambda_1$
then we have (\ref{rel-12}) and the desired implication 
 ``if the first process is hegemonic then the second process is also 
hegemonic" is correct. 

Because the rate function $L(A,B)$ is superadditive (in fact, additive) we can
immediately deduce part (i) of the Proposition (the constant $c$ is just
time-scaling) from the bond-percolation fact.

To prove (ii), fix an arbitrary reference point in the plane and consider the
perimeter ($X_t$, say) of the region $\RR_t$ containing the reference point at time
$t$.  The definition of the \EP\ gives growth dynamics as follows, where $(B_i)$ are
the regions adjacent to $\RR_t$, and $\FF_t$ denotes the history of the entire process up to time $t$.
\begin{eqnarray*}
E (dX_t | \FF_t)
&=& \sum_i (\peri(\RR_t \cup B_i) - \peri(\RR_t)) \lambda(\RR_t,B_i) \ dt \\ &\leq &
\sum_i \peri(B_i) \lambda(\RR_t,B_i) \ dt\\
& \leq & C \ \peri(\RR_t)\ dt  \mbox{ by assumption (ii)}\\
& = & C X_t\ dt .
\end{eqnarray*}
So $\frac{d}{dt} EX_t \leq C \ E X_t$ and so 
$EX_t \leq e^{Ct} EX_0$.  
So regions have finite mean size at all times.

\vspace{0.12in}
An explicit example where (ii) applies is 
\begin{equation}
\lambda(A,B) = \frac{L(A,B)}{\max(\peri(A), \peri(B))} . 
\label{ex-explicit}
\end{equation}
Because in this case 
\[ \sum_j \lambda(A,B_j)  \peri(B_j)
\leq \sum_j L(A,B_j) = \peri(A) . \]
Any merger rate $\lambda^\prime$ which is slower than this $\lambda$ in sense
(\ref{rel-0}) will also be non-hegemonic, because  case (ii) will apply.  But, to
repeat the warning at (\ref{heg12}), in
general if a merger rate is slower in sense (\ref{rel-0}) than a known non-hegemonic
merger rate, then we
cannot immediately deduce that it is also non-hegemonic.

\subsection{Notes on section \ref{sec-heg}}

\paragraph{Details of the definition of \EP.}
It is natural to assume the initial configuration is statistically
translation-invariant, so that it will remain so at all times.  
In the section \ref{sec-sims}  simulations we started with the square grid for
simplicity, though for certain theoretical analyses 
(section \ref{sec-constant}) it is more convenient to use the hexagonal lattice, and
as a continuum model it is perhaps conceptually  more natural to start with a
non-lattice configuration such as the Voronoi tessellation on random points.

\paragraph{Details of the definition of hegemonic.}
Property (\ref{def-heg}) is imprecise in two ways, which we mention here without
seeking a more precise analysis. 
First, we are presuming that the choice of initial configuration does not affect
(\ref{def-heg}).  
Second, we are presuming that (\ref{def-heg}) has probability zero or one.
Both presumptions are plausible by comparison with standard percolation results.

\paragraph{Limitations of the coupling methodology.}
One might hope that the coupling methodology used to prove (i) might be useful more
widely to prove hegemony under different assumptions.  But, by considering small
squares on either side of a boundary line, one sees
 that any superadditive rate function 
$\lambda$ must satisfy the hypothesis of (i).  Since we need
superadditivity to use the coupling, one cannot handle any essentially different
cases this way. 

\paragraph{Limitations of the methodology for (ii).}
The result in (ii) remains true if we replace ``perimeter" by ``area" or indeed any
subadditive functional, but we don't know any interesting cases that can be handled
by this variant and not by the stated form of (ii).

\paragraph{Why might the implication (\ref{heg12}) fail?}
This very sketchy outline indicates why (\ref{heg12}) is not obvious.
Consider a hegemonic \EP\ with rate $\lambda_1$, starting with the square grid.  
Write $\RR_t$ for the region containing a reference  point.  
Take a sequence $m_1 , m_2 , \ldots$ such that (with probability one) only finitely many values $m_j$ are seen as values of 
$\area(\RR_t)$.  Then we should be able to modify $\lambda_1$ by redefining $\lambda_1(A,B)$ to be small whenever
$\area(A)$ or $\area(B)$ equals some $m_j$, preserving the hegemonic property and the ``only finitely many values $m_j$ are seen" property.
Now take $\lambda_2$ as the modification of $\lambda_1$ in which rates $\lambda_2(A,B)$ are made very large when 
$\area(A) + \area(B) = $ some $m_j$.  This modification should force the sequence $\area(\RR_t)$ (for the $\lambda_2$ \EP) to re-enter the sequence 
$(m_j)$ infinitely often and (having slow merger rates therein) to be non-hegemonic.

\section{Analysis of the case $\lambda = 1$}
\label{sec-constant} 
In this section we present a non-rigorous, though convincing, analytic argument that
the constant rate case
($\lambda(A,B) = 1$ for adjacent $A,B$) is hegemonic, and this is
supported by
simulations -- see section \ref{sec-sims}.

Note that a special feature of this case is that any initial boundary line will
still be a piece of some boundary line at time $t$ with probability exactly
$e^{-t}$.  However, since these events for different initial lines are dependent, it
is not clear how to exploit this formula directly. 

It is convenient (the reason is explained below) to take the initial configuration
to be the hexagonal
lattice, scaled so that boundary edges have length $1$.  Consider a circuit $\CC$ which
encloses a reference point $b$.  Write 
$\RR_t$ for the boundary of the region of the \EP\ containing $b$ at time $t$.  As
with the standard Peierls contour  method \cite{gri99,MR2283880} in percolation, to
prove hegemonic it is enough to prove that 
\[
\int_0^\infty \sum_{\CC} P(\RR_t = \CC) \ dt  < \infty  \]
where the sum is over all circuits enclosing $b$.  
The length of a circuit must be even and $\geq 6$: write $2n$ for the length.  The
number of possible circuits enclosing $b$ of length $2n$ can be bounded as order
$2^{2n}$ without taking into account their
self-avoiding property; but taking this into account reduces the bound to order 
$2^{(2 - \delta)n}$ for some $\delta > 0$  (see \cite{MR2065578} for 
discussion of the value of $\delta$, which is not important for our calculation).
So if we can derive an upper bound $p_{2n}(t)$ for the probability that any particular contour
of length $2n$
is the boundary of $\RR_t$, then it is enough to show 
\begin{equation}
\sum_n 2^{(2 - \delta)n} \int_0^\infty p_{2n}(t) \ dt  < \infty .
\label{to-show}
\end{equation}

To obtain such an upper bound  $p_{2n}(t)$ consider Figure 2, which shows part of a circuit $\CC$ 
which is present at time $t$, and shows the regions whose boundaries include part of
$\CC$: shown are interior regions $A,B,C,D$ and exterior regions
$W,X,Y,Z$.

\setlength{\unitlength}{0.45in}
\begin{picture}(6,6)(-1,-4)
\put(0,0){\line(1,0){1}}
\put(1,0){\line(100,173){0.5}}
\put(1,0){\line(100,-173){0.25}}
\put(1.5,0.866){\line(1,0){1}}
\put(2.5,0.866){\line(100,173){0.25}}
\put(2.5,0.866){\line(100,-173){0.5}}
\put(3,0){\line(1,0){1}}
\put(4,0){\line(100,173){0.25}}
\put(4,0){\line(100,-173){0.5}}
\put(4.5,-0.866){\line(-100,-173){0.5}}
\put(4,-1.732){\line(100,-173){0.5}}
\put(4,-1.732){\line(-1,0){0.5}}
\put(4.5,-2.6){\line(-100,-173){0.25}}
\put(4.5,-2.6){\line(1,0){1}}
\put(5.5,-2.6){\line(100,173){0.25}}
\put(5.5,-2.6){\line(100,-173){0.5}}
\put(6,-3.466){\line(1,0){1}}
\put(-0.3,-0.05){$\CC$}
\put(7.1,-3.6){$\CC$}
\put(0.4,-0.4){A}
\put(1.2,1.1){W}
\put(2.5,-0.4){B}
\put(3.2,0.5){X}
\put(3.6,-2.5){C}
\put(4.8,-1.4){Y}
\put(5.3,-3.2){D}
\put(6.1,-2.8){Z}
\put(1,-3){b}
\end{picture}

{\bf Figure 2.}  {\small Part of a contour $\CC$ enclosing a point $b$, illustrating
regions adjacent to the contour.}

\vspace{0.16in}
\noindent
This picture can change in two ways at time increases. \\
(i) Two regions that are adjacent along the circuit -- either exterior regions such
as $X,Y$ or interior regions such as $B,C$ -- may merge into one.   This preserves
the circuit and (typically) decreases the number of regions adjacent to the circuit
by $1$. \\
(ii) An interior and exterior region -- such as $B,X$ -- may merge, destroying the
circuit.  Note this is where it is convenient to start with the hexagonal lattice --
the same  position on $\CC$ cannot separate two exterior regions and two interior
regions, so with $i$ exterior and $j$ interior regions, the number  of adjacent interior-exterior region pairs equals
$i+j$, provided neither $i$ nor $j$ equals $1$.

\vspace{0.1in}
\noindent
We continue the analysis assuming the ``typical" behaviour above always occurs --
 see Discussion below.
Because each possible merger occurs at rate $1$, the process 
$(I_t,J_t)$ where
\[ I_t = \mbox{ number of exterior regions adjacent to the circuit} \] \[ J_t =
\mbox{ number of interior regions adjacent to the circuit} \] behaves as the
continuous-time Markov chain with transition rates  
\begin{eqnarray*}
 (i,j) &\to& (i-1, j) : \mbox{ rate $i$}  \\ 
      &\to& (i, j-1) : \mbox{ rate $j$ } \\ 
     &\to& \mbox{ destroyed} : \mbox{ rate $i+j $} .
\end{eqnarray*}
These rates correspond to the three possibilities that adjacent exterior regions merge; or adjacent interior regions merge; or an exterior and interior region merge.  
The latter possibility is the case that the circuit no longer uses only boundary edges of the empire process, and so in particular  cannot be the boundary of the region 
containing $b$.

Recall we are studying a circuit of length $2n$ in the hexagonal lattice, and it is easy to check such a circuit in this initial lattice has 
$n+3$ exterior and $n-3$ interior regions, so the initial state of the Markov chain is 
$(I_0,J_0) = (n+3,n-3)$.  
Strictly, these rates only apply when $i,j \geq 2$, but allowing 
``fictitious" transitions $1 \to 0$ simplifies the analysis (see
Discussion below).  
We can compare the  Markov chain  $(I_t,J_t)$ to the chain $(I^*_t,J^*_t)$ defined
in the same way but without
the ``destroyed" possibility.  
Because in the process $(I_t,J_t)$ the next transition has
chance exactly $1/2$ to be ``destroyed", 
we find 
\begin{equation}
 P(I_t = i, J_t = j) = 2^{i+j - 2n} P(I^*_t = i, J^*_t = j) .\label{I*I}
\end{equation}
We are interested in $p_{2n}(t) =  P( J_t = 1 \mbox{ or } 0)$ 
(the $0$ captures the fictitious transition), and combining (\ref{I*I}) and
(\ref{to-show}) we see that we need to prove 
\begin{equation}
\sum_n 2^{ - \delta n}  \sum_{i = 1}^{n-3} 2^i \int_0^\infty\
P_{n+3,n-3}(I^*_t = i, J^*_t = j \mbox{ or } 0)
\ dt  < \infty  \label{to-show-2}
\end{equation}
where we wrote $P_{n+3,n-3}(\cdot)$ as a reminder that the chain 
$(I^*_t,J^*_t)$ starts at $(n+3,n-3)$.

Now if the chain $(I^*_t,J^*_t)$ ever hits state $(i,j)$ then it remains in that
state for mean time $1/(i+j)$, so we can write 
\[ \int_0^\infty\ P_{n+3,n-3}(I^*_t = i, J^*_t = j \mbox{ or } 0)
\ dt =  (i+1)^{-1} q_{n+3,n-3}(i,1) + i^{-1} q_{n+3,n-3}(i,0) \]
where $ q_{n+3,n-3}(i,j)$ is the chance that the embedded discrete-time jump chain 
$(\hat{X}_s, \hat{Y}_s)$ ever hits state $(i,j)$.  But this jump chain is simply the
chain obtained by drawing without replacement from 
a box which initially has $n+3$ balls labeled "exterior" and $n-3$ balls labeled
``interior", and writing $(\hat{X}_s, \hat{Y}_s)$ for the number of remaining balls
of each type after $s$ draws.  In particular, by considering reversed order of
draws 
\[ q_{n+3,n-3}(i,1) = (i+1) \ \frac{n-3}{2n} \ \frac{(n+3)_i}{(2n-1)_i} \leq C (i+1)
2^{-i} \]
for some constant $C$, with a smaller similar bound for
$q_{n+3,n-3}(i,0)$.  
We now see the inner sum in (\ref{to-show-2}) is bounded by order $n$ and thus the outer
sum is indeed finite.

\subsection{Discussion of approximations above.}
What is assumed in setting up the rates for the continuous-time Markov chain is that
the $j$ successive exterior (and similarly for interior) regions touching the
circuit $\CC$ \\
(a) are distinct \\
(b) meet only if they are consecutive pairs along $\CC$ \\
(c) $j \geq 3$. \\
Under these assumptions, the combined merger rate is indeed $j$.   If instead of (c)
we have $j = 2$ then the merger rate is $1$, and if $j = 1$ the rate is $0$; however
it is not hard to see that these changed rates for $j = 1,2$ don't affect the
``finiteness" result (\ref{to-show}).

But the other simplifying assumption are more significant.  If in Figure 2 the regions $W$ and $Y$ are
adjacent, by meeting behind region $X$, then the combined merger rate is larger than
$j$; moreover if they do merge then (a) fails for the merged region.  
The resulting ``combinatorial explosion" of possibilities for successive regions along the boundary being distinct or identical 
seems very hard to analyze rigorously.

\section{Simulation results}
\label{sec-sims}
We will show simulation results for the following four rate functions.

\noindent
{\bf Model 1:}  $\lambda(A,B) = 1 $.

\noindent
{\bf Model 2:}  $\lambda(A,B) = \area(A) \times \area(B) $.

\noindent
{\bf Model 3:}  $\lambda(A,B) = L(A,B) $.

\noindent
{\bf Model 4:}  $\lambda(A,B) = 1/(\area(A) \times \area(B)) $.

\noindent
The quantities we calculated were the spatial averages
\begin{eqnarray*}
S(t) &=& \mbox{ average (area of empire)}^2\mbox{ per unit area } \\
D(t) &=& \mbox{ average number of empires per unit area. } 
\end{eqnarray*}
Note the alternative interpretations:
 the average area of an empire (for uniform random choice of empire) equals $1/D(t)$, 
 while the average area of the empire containing a random point (i.e. the average when empires are chosen wih probability proportional to area) 
 equals $S(t)$.
 We started at $t = 0$ with the square grid, 
 so $S(0) = D(0) =1$; as $t$ increases, $S(t)$ increases and $D(t)$ decreases.  
 In order to compare different models it is convenient to plot the curve $(D(t),S(t))$, regarded as a function $S(D)$.
 This is a natural way to summarize the non-hegemonic phase.
 If the model is hegemonic then, as $D$ decreases from $1$,  $S(D)$ will increase to infinity at some critical value $\Dcrit > 0$, whereas for a  non-hegemonic model 
 $\Dcrit = 0$.
 
 There is a different way to summarize the hegemonic phase.  Restrict the model to a large finite region of area $A$.  
 Then as $D$ decreases from $1$ to $0$, the quantity $S(D)/A$  increases from $1/A$ to $1$, and $S(D)/A$  
 equals the chance that two random points in the region are in the same empire.  
 In the $A \to \infty$ limit, the function $S(D)/A$ will become a function $f(D)$ analogous to the {\em percolation function} in percolation theory, with 
 $f(D) = 0$ for $D > \Dcrit$ and $f(D) > 0$ for  $D < \Dcrit$.  
 
These two ways to think about $\Dcrit$ are shown in Figure 3 for the four models.  
 The upper graph in each pair shows $S(D)$, and the lower graph shows $S(D)/A$, for two models.  
 The points $\Dcrit$ are where (in the $A \to \infty$ limit) the upper graph goes to infinity and where the lower graph leaves zero. 
 Data are from simulations on an $81 \times 81$ grid.
 
 \newpage
      \vspace*{-1in}
   \includegraphics[scale=0.74]{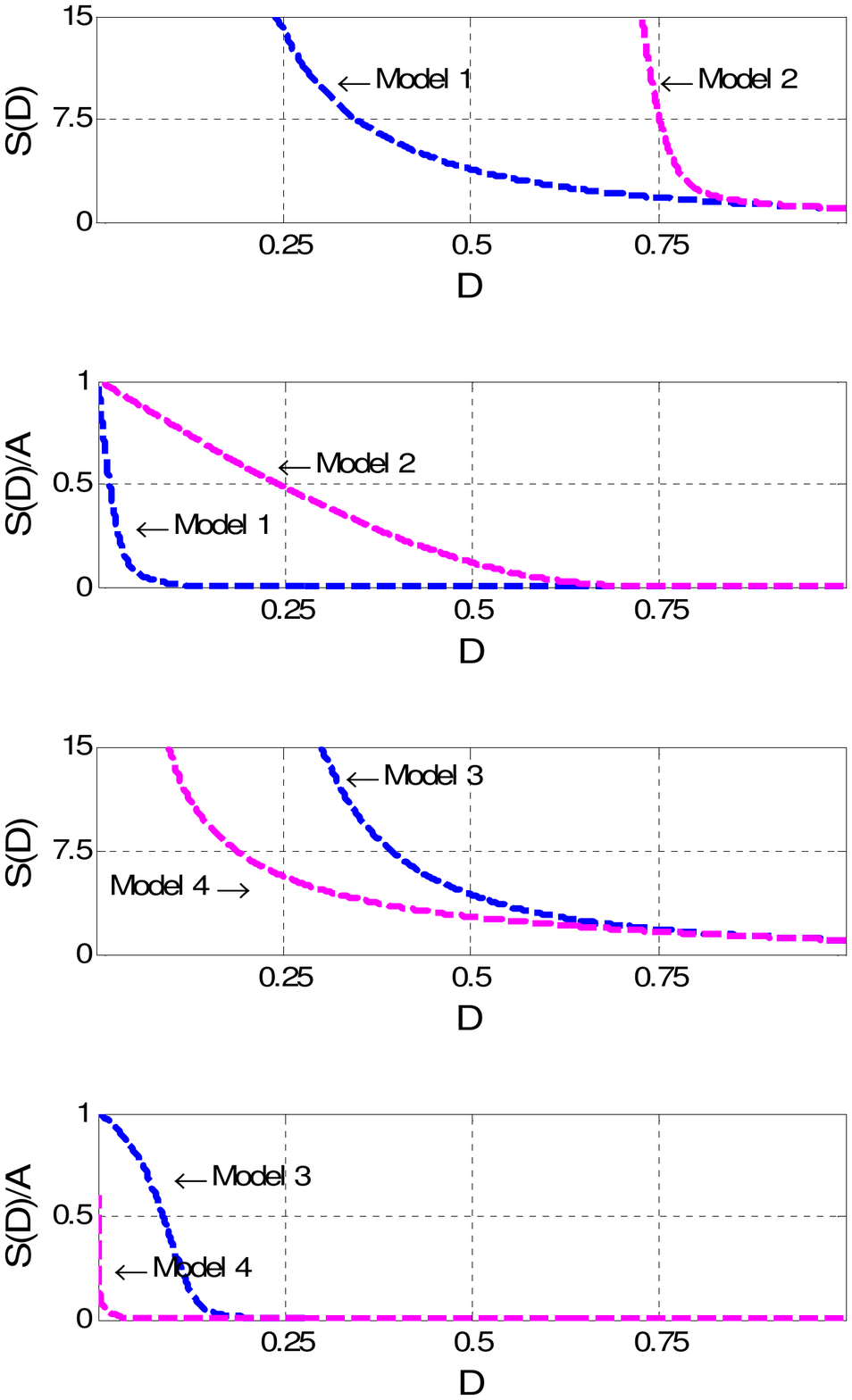}

      \vspace*{-1in}
 {\bf Figure 3.} {\small Simulation results for 4 models, explained in text.}
 
   \newpage
     \vspace*{-2in}
   \includegraphics[scale=0.64]{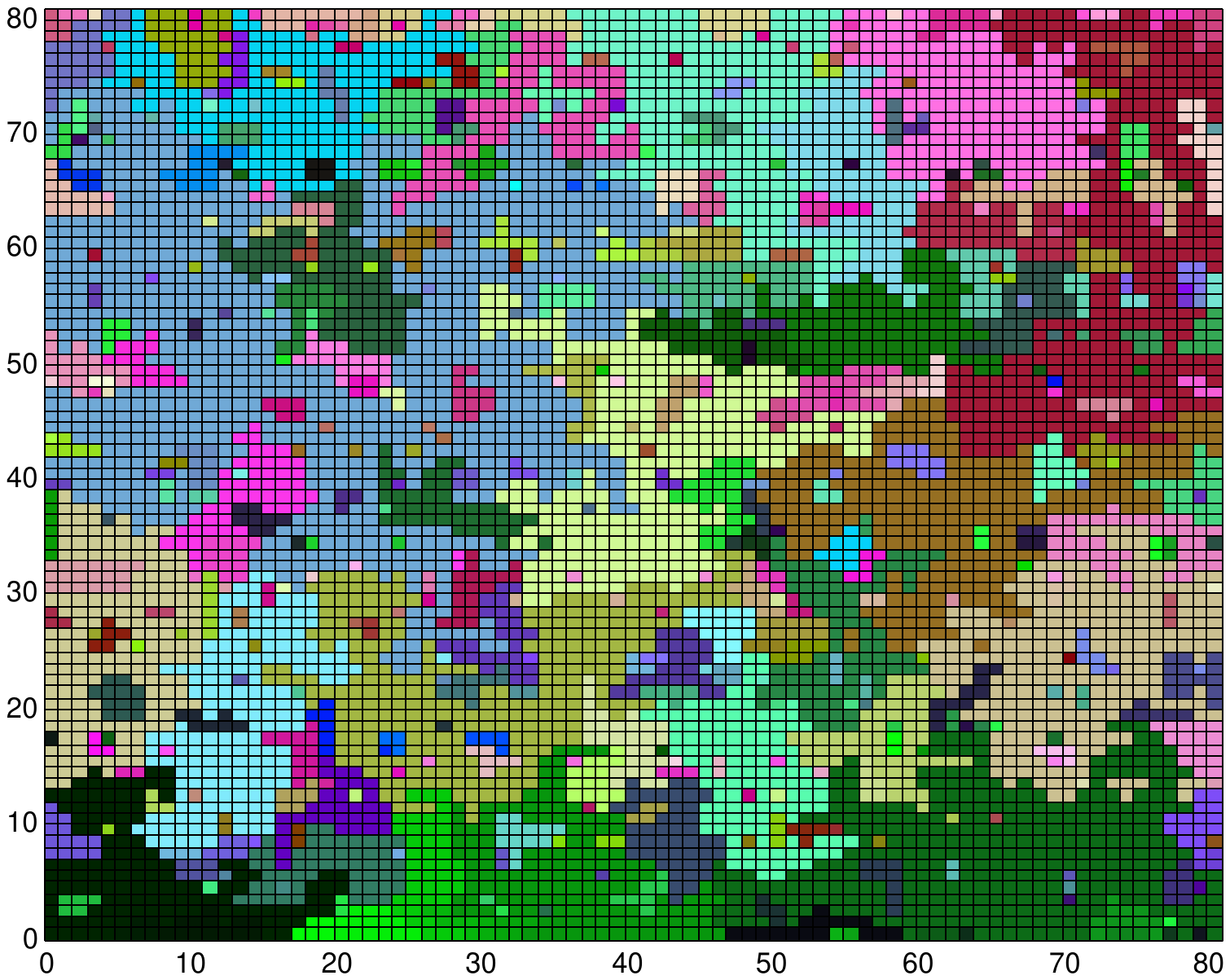}
   
   \vspace*{-3.4in}

     \includegraphics[scale=0.64]{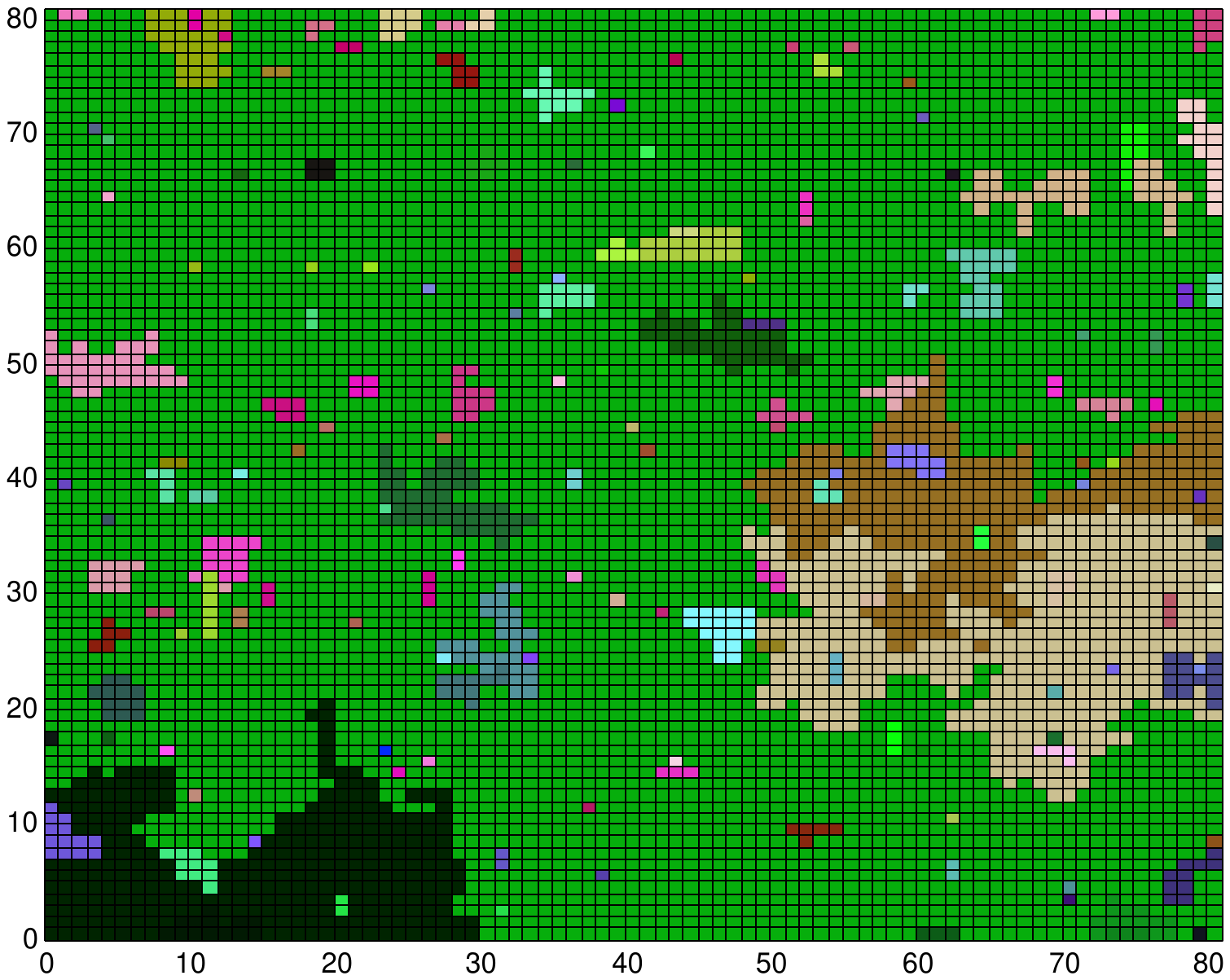}
     
       \vspace*{-1.7in}
     {\bf Figure 4.}  
     {\small The emergence of a giant component, in the case $\lambda = 1$.  
     The top panel shows a configuration at  $D = 0.075$ and the bottom at $D = 0.025$.
     }

 \newpage
The data from Figure 3 suggests the following rough numerical estimates, which we then compare to theory.  
The scaling exponents $\gamma$ will be discussed below.
   
\noindent
{\bf Model 1: [$\gamma = 0$]:}  $\Dcrit \approx  0.15$.  
This is consistent with the non-rigorous argument for hegemony ($\Dcrit >   0$). 

\noindent
{\bf Model 2: [$\gamma = 2$]:}  $\Dcrit \approx  0.6$.  In this case Proposition \ref{P1}(a) applies, so we know $\Dcrit >   0$. 

\noindent
{\bf Model 3:  [$\gamma = 1/2$]:}  $\Dcrit \approx  0.2$.  
This is the ``bond percolation" case, so we know $\Dcrit  > 0$.

\noindent
{\bf Model 4:  [$\gamma = -2$]:}  $\Dcrit \approx  0$.  In this case Proposition \ref{P1}(b) applies, so we know $\Dcrit =   0$. 

\vspace{0.1in}
\noindent
Figure 4 shows emergence of the giant component in the case $\lambda = 1$.   
The visual appearance is rather different from that of near-critical bond percolation, in that the large components appear less ``fractal" and that  distinct moderately large components coexist over longer time periods.

\vspace*{-0.2in}

\paragraph{Scaling exponents.}
In the context of mean-field coalescence, a kernel such that $K(cx,cy) = c^\gamma K(x,y)$ is said to have 
 scaling exponent $\gamma$. 
It has long been understood, mostly non-rigorously (but see \cite{MR2217655} for references to recent rigorous work) that 
the coalescence process should be non-gelling if $\gamma \leq 1$ but gelling if $\gamma > 1$.

We can define a scaling exponent $\gamma$ analogously for the \EP:
$\lambda(cA,cB) = c^{2\gamma} \lambda(A,B)$, where $cA$ denotes linear scaling by a factor $c$.  
The data above is consistent with the possibility that there is some critical value $\gamcrit$ 
such that \EP es are typically hegemonic for $\gamma > \gamcrit$ and are typically non-hegemonic for $\gamma < \gamcrit$.  
However, the rate function at (\ref{ex-explicit}) which is known to be non-hegemonic has $\gamma = 0$, as does the $\lambda = 1$ case presumed to be hegemonic,
so it may be that scaling exponents are not so definitive for \EP es.

\end{document}